\documentclass[runningheads]{llncs}
\usepackage{amsfonts,amssymb,latexsym}
\def\MW{Mordell-\mbox{\kern-.16em}Weil}
\def\NeS{N\'eron-Severi}
\newcommand{\nl}{\hspace*{\fill}\\}
\newcommand{\be}{\begin{equation}}
\newcommand{\ee}{\end{equation}}
\newcommand{\bea}{\begin{eqnarray}}
\newcommand{\eea}{\end{eqnarray}}
\newcommand{\ra}{\rightarrow}
\newcommand{\Ra}{\Rightarrow}
\newcommand{\IFF}{\underline{iff}}

\newcommand{\C}{\mathbf{C}}
\newcommand{\G}{\mathbf{G}}
\newcommand{\Q}{\mathbf{Q}}

\newcommand{\Z}{\mathbf{Z}}
\newcommand{\PP}{\mathbf{P}}
\newcommand{\EE}{\mathcal{E}}
\newcommand{\dd}{\mathbf{d}}
\newcommand{\LL}{\mathcal{L}}
\newcommand{\NN}{\mathbf{N}}
\newcommand{\lala}{{\mbox{\boldmath$\lambda$}}}
\newcommand{\XX}{\mathcal{X}}
\newcommand{\Ht}{\hat{h}}
\newcommand{\Htmin}{\Ht_{\min}}
\newcommand{\naive}{na\"\i ve}
\newcommand{\End}{\mathop{\rm End}\nolimits}

\newcommand{\Jac}{\mathop{\rm Jac}\nolimits}
\newcommand{\rk}{\mathop{\rm rk}\nolimits}
\newcommand{\NS}{\mathop{\rm NS}\nolimits}
\newcommand{\0}{^{\phantom0}}
\newcommand{\9}{_{\phantom0}}
\newcommand{\z}{{\phantom0}}
\newcommand{\Qed}{\nopagebreak[4]\hfill
\rule{2mm}{2.5mm} \bigskip \pagebreak[2]}
\begin{document}

\title{Points~of~low~height~on~elliptic~curves~and~surfaces\\
I: Elliptic surfaces over~$\PP^1$ with small~$d$}
\titlerunning{Points of low height on elliptic curves and surfaces I}
\author{Noam D.~Elkies}
\institute{Department of Mathematics, Harvard University, Cambridge, MA
02138 USA}
\maketitle

\begin{abstract}
For each of $n=1,2,3$ we find the minimal height $\Ht(P)$ of
a nontorsion point~$P$\/
of an elliptic curve~$E$\/ over $\C(T)$ of discriminant degree $d=12n$
(equivalently, of arithmetic genus~$n$),
and exhibit all $(E,P)$ attaining this minimum.
The minimal $\Ht(P)$ was known to equal $1/30$ for $n=1$
(Oguiso-Shioda) and $11/420$ for $n=2$ (Nishiyama),
but the formulas for the general $(E,P)$ were not known,
nor was the fact that these are also the minima for an elliptic curve
of discriminant degree~$12n$ over a function field of any genus.
For $n=3$ both the minimal height ($23/840$)
and the explicit curves are new.  These $(E,P)$
also have the property that that $mP$\/ is an integral point
(a point of na\"\i ve height zero) for each $m=1,2,\ldots,M$,
where $M=6,8,9$ for $n=1,2,3$; this, too, is maximal in each
of the three cases.
\end{abstract}

\vspace*{4ex}

{\large\textbf{1.~Introduction.}}

\textbf{1.1 Statement of results.}
Let $K$\/ be a function field of a curve~$C$\/ of genus~$g$
over a field~$k$ of characteristic zero,\footnote{
  One can also usefully define the canonical height etc.\
  in positive characteristic, but we need to use the ABC conjecture
  for~$K$\/ and thus must assume that $K$\/ has characteristic zero.
  }
and $E$\/ a nonconstant elliptic curve over~$K$.
Let $d$\/ be the degree of the discriminant of~$E$\/
(considered as a divisor on~$C$\/),
a natural measure of the complexity of~$E$\/;
and let $\Ht: E(K) \ra \Q$ be the canonical height.
Necessarily $12|d$\/; in fact it is known that $d=12n$
where $n$ is the arithmetic genus of the elliptic surface~$\EE$
associated with~$E$.  It is not hard to show that, given~$d$,
the set of numbers~$H$\/ that can occur as the canonical height
of a rational point on~$E$\/ is discrete.
In particular, for each $d=12n$ there is a minimal positive height
$\Htmin(d)$, and also a minimal positive height $\Htmin(g,d)$
for elliptic curves over function fields of genus~$g$
(except for $g=d=0$, when $E$\/ is a constant curve over~$\PP^1$
and thus has no points of positive height).
It is thus a natural problem to compute or estimate these numbers
$\Htmin(d)$ and $\Htmin(g,d)$.  This paper is the first of a series
concerned with different aspects of this problem.

In this paper we determine $\Htmin(12n)$ for $n=1,2$
and $\Htmin(0,12n)$ for $n=1,2,3$.
Since we are working in characteristic zero, we may assume $k=\C$,
when every \hbox{genus-zero} curve is isomorphic to $\PP^1$
and its function field is isomorphic to $\C(T)$.

\begin{theorem}
\label{thm:n=1}
i) (Oguiso-Shioda \cite{OgSh}) $\Htmin(0,12)=1/30$. \nl
ii) $\Htmin(12)=1/30$.
Moreover, let $E$\/ be an elliptic curve with $d=12$
over a complex function field~$K$, and $P\in E(K)$.
Then the following are equivalent:
(a)~$\Ht(P)=1/30$;
(b)~Each of $P,2P,3P,4P,5P,6P$\/ is an integral point on~$E$;
(c)~$K \cong \C(T)$, and $(E,P)$ is equivalent to the curve
\be
\label{eq:E1}
E_1(q) :
Y^2 + (s'-(q+1)s) XY + qss'(s-s') Y = X^3 - qss' X^2
\ee
over the $(s:s')$ line with the rational point $P : (X,Y) = (0,0)$,
for some $q\in\C$ other than $0$ or~$1$.
\end{theorem}

\begin{theorem}
\label{thm:n=2}
i) (Nishiyama \cite{Nishiyama}) $\Htmin(0,24)=11/420$. \nl
ii) $\Htmin(24)=11/420$.
Moreover, let $E$\/ be an elliptic curve with $d=24$
over a complex function field~$K$, and $P\in E(K)$.
Then the following are equivalent:
(a)~$\Ht(P)=11/420$;
(b)~$mP$ is an integral point on~$E$\/ for each $m=1,2,\ldots,8$;
(c)~$K \cong \C(T)$, and $(E,P)$ is equivalent to the curve
\bea
& \!\!\!\!\!\!\!
E_2(u) :
&
Y^2 + (r^2-{r'}^2+(u-2)rr') XY
\nonumber \\
&&
- \; r^2 r' (r+r') (r+ur') (r+(u-1)r') Y
\label{eq:E2} \\
&&
= X^3 - r r' (r+r') (r+ur') X^2
\nonumber
\eea
over the $(r:r')$ line with the rational point $P : (X,Y) = (0,0)$,
for some $u\in\C$ other than $0,1$.
\end{theorem}

\begin{theorem}
\label{thm:n=3}
i) $\Htmin(0,36)=23/840$. \nl
ii) Let $E/\C(T)$ be an elliptic curve with $d=36$,
and $P$\/ a rational point on~$E$.  Then the following are equivalent:
(a)~$\Ht(P)=23/840$;
(b)~$mP$ is an integral point on~$E$\/ for each $m=1,2,\ldots,9$;
(c)~$(E,P)$ is equivalent to the curve
\bea
& \!\!\!\!\!\!\!
E_3(A) :
&
Y^2 + (At^3 + (1-2A) t^2 t' - (A+1) t{t'}^2 - {t'}^3) XY
\nonumber \\
&&
- \; t^3 t' (t+t') (At+t') (At+(1-A)t') (At^2+tt'+{t'}^2) Y
\label{eq:E3} \\
&&
= X^3 - t t' (t+t') (At+t') (At^2+tt'+{t'}^2) Y
\nonumber
\eea
over the $(t:t')$ line with the rational point $P : (X,Y) = (0,0)$,
for some $A\in\C$ other than $0,1$.
\end{theorem}

The values of $\Htmin(12)$ and $\Htmin(24)$ are new.
Note that we do not claim to determine $\Htmin(36)$.
As indicated, the values of $\Htmin(0,12)$ and $\Htmin(0,24)$
(the first parts of Theorems~\ref{thm:n=1} and~\ref{thm:n=2})
were already known, but were obtained using techniques
that are specific to the geometry of rational and K3 elliptic surfaces
and do not readily generalize past $n=2$.
Our approach lets us treat all three cases uniformly,
and in principle lets us determine $\Htmin(0,12n)$ for any~$n$, though
the computations rapidly become infeasible as $n$ grows beyond~$3$.
The minimizing $(E,P)$ had not been previously exhibited,
except for a single case of a rational elliptic surface
with a section of height~$1/30$ obtained by Shioda in a later paper
\cite{Sh:exist}, which we will identify with $E_1(4/5)$.

The connections with integral multiples of~$P$\/
(see statement~(b) of part~(ii) of each Theorem) are also new.
We do not expect them to persist past \hbox{$n=3$,}
and in fact find that for $n=4$ the largest number
of consecutive integral multiples occurs for $(E,P)$ with
$\Ht(P) = 19/630$ or $13/360$, whereas
$\Htmin(0,48) \leq 41/1540 < 19/630 < 13/360$.
We shall say more about integrality later;
for now we content ourselves with the following remarks.
A point on an elliptic curve over a function field~$k(C)$ is said to be
integral if it is a nonzero point whose \naive\ height vanishes.
Geometrically, if we regard $E$\/ as an elliptic surface~$\EE$ over~$C$,
and a rational point $P\in E(K)$ as a section~$s_P$ of~$\EE$,
this means that $s_P$ is disjoint from the zero-section $s_0$ of~$\EE$.
Since $g=0$ in our case, we can give an explicit algebraic
characterization of integrality.  Write $E$\/ in extended
Weierstrass form as
\be
Y^2 + a_1 XY + a_3 Y = X^3 + a_2 X^2 + a_4 X + a_6
\label{eq:Weier}
\ee
where each $a_i$ is a homogeneous polynomial of degree $i \cdot n$
in two variables.  Then a rational point $(X,Y)$ is integral
if $X,Y$\/ are homogeneous polynomials of degrees $2n,3n$ respectively.
The equation (\ref{eq:Weier}) depends on the choice of coordinates
$X,Y$\/ on~$E$\/; replacing $X,Y$\/ by
\be
\label{eq:ellchange}
\delta^2(X + \alpha_2), \qquad
\delta^3(Y + \alpha_1 X + \alpha_3)
\ee
(some $\alpha_i$ and nonzero~$\delta$) yields an isomorphic curve.
If moreover $\delta\in\C^*$ and each $\alpha_i$ is a homogeneous
polynomial of degree $i \cdot n$ then the new equation for~$E$\/
has the same discriminant degree and the same integral points.

\textbf{1.2 Outline of this paper.}
For each $n=1,2,3$ we prove Theorem~$n$,
except for the implications (a),(b)$\Ra$(c) of part~(ii),
which require different methods that we defer to a later paper.
Our proofs use the following ingredients:
\begin{itemize}
\item $\Ht(mP) = m^2 \Ht(P)$ for all $m\in\Z$.
\item If $mP \neq 0$ then
\be
\label{H=h+}
\Ht(mP) = h(mP) + \sum_v \lambda_v(mP),
\ee
where $h(\cdot)$ is the \naive\ height and the sum extends over
all places $v\in C(\C)$ lying under singular fibers~$E_v$ of~$E$.
(All places of~$K$\/ are of degree~$1$ thanks to our use of the
algebraically closed field~$\C$ for~$k$.)
The local corrections $\lambda_v(mP)$ are described further below.
\item The \naive\ height takes values in $\{0,2,4,6,\ldots\}$,
and satisfies $h(m'P) \leq h(mP)$ for any integers $m,m'$
such that $m'|m$ and $mP \neq 0$.
\item Each local correction $\lambda_v(mP)$ depends only
on the Kodaira type of the fiber~$E_v$ and on the component
of $E_v$ meeting~$P$.  We shall call this component $c_v$.
The values of $\lambda_v(\cdot)$ are known explicitly
for all Kodaira types and each possible component,
see for instance \cite[Thm.~5.2]{Sil:ht}.
\item Finally, the condition that $E$\/ have discriminant degree $d=12n$
imposes two conditions on the Kodaira types of the singular fibers.
The first condition is
\be
\label{d}
d = \sum_v d_v,
\ee
where $d_v$ is the local discriminant degree of~$E_v$.
This allows only finitely many collections of fiber types.
The second condition follows from an inequality
due to Shioda \cite[Cor.~2.7 (p.30)]{Sh:Ineq},
and eliminates some of these collections that have too few fibers.
According to this condition, if a nonconstant elliptic curve
of discriminant degree~$d$\/ over a function field $K=\C(C)$
has a nontorsion point then the conductor degree of the curve
strictly exceeds $(d/6) + \chi(C)$.
Here $\chi(C)=2-2g$ is the Euler characteristic of~$C$.
The conductor degree may be defined as the number of multiplicative
fibers plus twice the number of additive fibers; thus it is also
a sum of invariants of the singular fibers.
When $(g,d)=(0,12n)$ we have $\chi(C)=2$ and $d/6=2n$,
so the conductor degree is at least $2n+3$.
\end{itemize}
We shall refer to these constraints as the ``combinatorial conditions''
on $\Ht(P)$, $h(mP)$, and the collection of $(E_v,c_v)$
that arise for $(E,P)$.  (For other uses of such conditions
to obtain lower bounds on heights, see for instance
\cite{HS,Sil:Lehmer} and work referenced in these sources.)
In general the combinatorial conditions yield
only a lower bound on $\Htmin(0,12n)$, because they allow
some possibilities that do not actually occur for any $(E,P)$.
But for each of $n=1$, $2$, and $3$ this lower bound turns out
to be attained by some $(E,P)$ over~$\C(T)$,
namely those exhibited in statement~(c) of part (ii) of Theorem~$n$.
(Note that we do not yet need to derive the formulas for these $(E,P)$,
nor to prove that they are the only ones possible.)
Moreover, using (\ref{H=h+}) we can check that $\Ht(P)=\Htmin(0,12n)$
if and only if the \naive\ height $h(mP)$ vanishes
for all $m$ up~to $6$, $8$, or~$9$ respectively.

Still, already at $n=1$ we see some redundancy.
The combinatorial conditions allow $\Ht(P)=1/30$
to be attained in any of five ways,
four of which are realized by the curves $E_1(q)$
of Theorem~\ref{thm:n=1} for suitable choices of~$q$.
Shioda's $E_1(4/5)$ has singular fibers of types
I$_5$, I$_3$, I$_2$, and~II.
(We specify the components $c_v$ later in the paper.)
The fibers of $E_1(-1)$ have types I$_5$, IV, I$_2$, and~I$_1$,
while those of $E_1(4)$ have types I$_5$, I$_3$, III, and~I$_1$.
In all other cases, the fibers of $E_1(q)$ have types
I$_5$, I$_3$, I$_2$, I$_1$, I$_1$: the first three at
$s=0$, $s'=0$, $s'=s$, and the last two at the roots of the quadratic
$
(q+1)^3 s^2 = (11q^2-14q+2) ss' + (q-1) {s'}^2.
$
When $q=4/5$, these roots coincide and the two I$_1$ fibers merge
to form a~II; likewise at $q=-1$ or $q=4$, one of the I$_1$ fibers
merges with the I$_3$ or I$_2$ fiber to form a IV or~III respectively.
(The one merger that does not occur is ${\rm I}_1 + {\rm I}_1 \ra {\rm I}_2$.)
But none of these degenerations changes $\Ht(P)$, nor any $h(mP)$,
nor the conductor degree~$N$.  In fact a fiber of type II, III, or~IV
contributes as much to our formulas for $\Ht(P), h(mP), N$\/
as a pair of fibers of types I$_1$ and I$_\nu$ ($\nu=1$, $2$, or~$3$).
Thus it is enough to minimize $\Ht(P)$ under the further assumption
that no fibers of type II, III, or~IV occur.
We find similar replacements for all components of
fibers of the remaining additive types
I$^*_\nu$, II$^*$, III$^*$, IV$^*$.  See Proposition~2.  This simplifies
the computation of the combinatorial lower bound on $\Htmin(0,12n)$:
instead of an exhaustive search over all combinations of $(E_v,c_v)$,
we need only try those for which each $E_v$ is multiplicative 
(of type $I_\nu$ for $\nu=d_v$).

We programmed the search over all partitions $\{d_v\}$ of~$12n$
in~{\sc gp}~\cite{Pari} and ran it on a Sun Ultra~60.
This took only a fraction of a second for $n=1$,
five seconds for $n=2$, and five minutes for $n=3$.
It took about an hour to carry out the same computation for $n=4$,
and about $20$ hours for $n=5$;
but the resulting bounds are probably not attained:
as we shall see in a later paper,
the required $(E_v,c_v)$ data impose more conditions
than the number of parameters needed to specify $(E,P)$.
We do produce explicit $(E,P)$ that show $\Htmin(0,48) \leq 41/1540$
and $\Htmin(0,60) \leq 261/10010$, and conjecture that these are
the correct values of $\Htmin(0,12n)$ for $n=4,5$.
We have not attempted to extend the computation past $n=5$.

\textbf{1.3 Coming attractions.}
Happily, the computation of the surfaces
(\ref{eq:E1},\ref{eq:E2},\ref{eq:E3})
not only completes the proofs of Theorems~1 through~3
but also points the way to further results and connections.
We outline these here, and defer detailed treatment
to a later paper in this series.
In each step of the computation we in effect
obtain a new birational model for the moduli space, call it~$\XX$,
of pairs $(E,P)$ consisting of an elliptic curve and a point on~it.
Our new parametrizations of this rational surface~$\XX$\/
have several other applications.  One is a geometric interpretation
of Tate's method for exhibiting the generic elliptic curve
with an \hbox{$N$-torsion} point: we readily locate the modular curves
X$_1(N)$ ($N \leq 16$) on~$\XX$, together with nonconstant
rational functions of minimal degree that realize each X$_1(N)$
as an algebraic curve of genus $\leq 2$.
Arithmetically, we can use our parametrizations of~$\XX$\/
to find $(E,P)$ over~$\Q$ (or over some other global field)
such that $P$\/ is a nontorsion point with small $\Ht(P)$,
and/or with many integral multiples in the minimal model of~$E$.
For instance, we prove that there are infinitely many $(E,P)/\Q$
such that $mP$\/ is integral for each $m=1,2,\ldots,11,12$.
Our numerical results for a isolated curves $(E,P)$ over~$\Q$
may be found on the Web at
{\sf http://www.math.harvard.edu/$\sim$elkies/low\_height.html}~.
They include new records for consecutive integral multiples
and for the Lang ratio $\Ht(P)/\log|\Delta_E|$.
We have $mP$\/ integral for each $m=1,2,\ldots,13,14$ for
\be
\label{eq:E14}
E: Y^2 + XY = X^3 - 139761580 X + 1587303040400,
\ee
an elliptic curve of conductor
$1029210 = 2 \cdot 3 \cdot 5 \cdot 7 \cdot 13^2 \cdot 29$,
and $P$\/ the nontorsion point $(X,Y) = (11480, 1217300)$;
and we find the curve
\be
\label{eq:EPmin}
Y^2 + XY = X^3 - 161020013035359930 X + 24869250624742069048641252
\ee
of conductor $3476880330 =
2 \cdot 3 \cdot 5 \cdot 7 \cdot 11 \cdot 23 \cdot 31 \cdot 2111$
with the nontorsion point $(-296994156, 6818852697078)$
of canonical height\footnote{
  There are two standard normalizations, differing by a factor of~$2$,
  for the canonical height of a point on an elliptic curve over~$\Q$.
  We use the larger one, which is the one
  consistent with our formulas for function fields.
  }
$\Ht(P) = .0190117\ldots < 1.691732 \cdot 10^{-4} \log|\Delta_E|$.
The curves (\ref{eq:E14},\ref{eq:EPmin}) are the specializations
of our formula~(\ref{eq:E3}) with
$(A, t/t') = (35/32, -8/15)$, $(33/23, 115/77)$.

Our simplified formula for $\Ht(mP)$ (Proposition~2)
also bears on the asymptotic behavior of $\Htmin(g,12n)$
for fixed~$g$ as $n\ra\infty$.
Hindry and Silverman~\cite{HS} used the combinatorial conditions
(except for the condition: $h(m'P) \leq h(mP)$ if $m'|m$)
to show that there exists $C>0$ such that
\be
\Ht(g,12n) \geq Cn - O_g(1),
\label{HS}
\ee
This proved the function-field case
of a conjecture of Lang~\cite[p.92]{Lang}.
The error terms $O_g(1)$ are effectively computed,
and can be omitted entirely if $g \leq 1$.
Hindry and Silverman also produce an explicit constant~$C$,
but it is quite small: about $7\cdot 10^{-10}$.
Their approach requires a point meeting every additive fiber
in its identity component,
which they achieved by working with $12P$\/ instead of~$P$,
at the cost of a factor of $1/12^2$ in~$C$.
Our results here let one apply the same methods directly to~$P$,
thus saving a factor of~$12^2$ and raising $C$\/ to about $10^{-7}$.
In a later paper
we show how to gain another factor of approximately $5000$,
raising the lower bound on $\liminf_n \Ht(g,12n)/n$ to $1/2111$.
This is within an order of magnitude of the correct value:
for all $n \equiv 0 \bmod 5$ we obtain $\Htmin(0,12n)\leq 261n/50050$
via base change from our $n=5$ example.

\vspace*{4ex}

{\large\textbf{2.~The \naive\ and canonical heights.}}

We collect here the facts we shall use about elliptic curves~$E$\/
over function fields~$K$\/ in characteristic zero, the associated
elliptic surface~$\EE$, and the \naive\ and canonical height functions
on~$E(K)$.

\textbf{2.1 The \naive\ height.}
The \textit{\naive\ height} $h(P)$ of a nonzero $P\in E(K)$
can be defined using intersection theory on the elliptic surface~$\EE$\/
associated to some model of~$E$.
Let $s_0$ be the zero-section of the elliptic fibration $\EE \ra C$,
and $s_P$ the section corresponding to~$P$.  Then
  $
  h(P) := 2 s_P \cdot s_0.
  $
Since we assumed that $P \neq 0$,
the sections $s_0,s_P$ are distinct curves on~$\EE$.
Hence their intersection number $s_P\cdot s_0$ is a nonnegative integer,
and $h(P)$ is a nonnegative even integer.
Moreover $h(P)=0$ if and only if $s_P$ is disjoint from $s_0$,
in which case we say that $P$\/ is an \textit{integral point}\/ on~$E$.

When $C=\PP^1$, we can give an equivalent algebraic definition of~$h(P)$
in terms of a Weierstrass equation of~$E$.  This definition emphasizes
the analogy with the canonical height in the more familiar case
of an elliptic curve over~$\Q$.  Recall that each coefficient~$a_i$
in the Weierstrass equation~(\ref{eq:Weier}) is a homogeneous polynomial
of degree $i \cdot n$ in the projective coordinates on~$\PP^1$.
Then the coordinates $x,y$ of a nonzero $P\in E(K)$
are homogeneous rational functions of degrees $2n,3n$.
If $x,y$ are written as fractions ``in lowest terms''$\!$,
as quotients of coprime homogeneous polynomials, then the denominators are
(up to scalar multiple) the square and cube of some polynomial~$\zeta$.
The roots of~$\zeta$, with multiplicity, are the images on~$\PP^1$
of the intersection points of $s_0$ and $s_P$.  Hence
  $
  s_P \cdot s_0 = \deg\zeta.
  $
Therefore $h(P)$ is the degree of the denominator~$\zeta^2$ of~$x$,
which is also the number of poles of~$x$ counted with multiplicity.
An integral point is one for which $\zeta$ is a nonzero scalar
and thus $x,y$ are homogeneous polynomials of degrees~$2n,3n$.

For an arbitrary base curve~$C$, the coefficients~$a_i$ are
global sections of $\LL^{\otimes i}$ for some line bundle~$\LL$ on~$C$,
and $x,y$ are meromorphic sections of $\LL^{\otimes 2},\LL^{\otimes 3}$.
The pole divisors of~$x,y$ are $2Z,3Z$
for some effective divisor~$Z$\/ on~$C$, whose degree is $s_P\cdot s_0$;
thus again $h(P)$ is the degree of the pole divisor~$2Z$\/ of~$x$,
and $P$\/ is integral \IFF\ $Z=0$
\IFF\ $x,y$ are global sections of $\LL^{\otimes 2},\LL^{\otimes 3}$.
A linear change of coordinates according to (\ref{eq:ellchange})
yields the same notion of integrality if and only if $\delta\in\C^*$
and $\alpha_i \in \Gamma(\LL^{\otimes i})$ for each~$i$.

We shall need one more property of the \naive\ height beyond
its relation with the canonical height and the fact that
$h(mP) \in \{0,2,4,6,\ldots\}$ ($mP\neq0$):

\begin{lemma}
\label{lemma:ht_ineq}
Let $P$\/ be a point on an elliptic curve over $k(C)$,
and let $m,m'$ be any integers such that $m'|m$ and $mP \neq 0$.
Then $h(m'P) \leq h(mP)$.
\end{lemma}

\textit{Proof}\/:
Each point of $s_{m'P} \cap s_0$ is also a point of intersection
of $s_{mP}$ with $s_0$, to at least the same multiplicity.
Hence $s_{m'P} \cdot s_0 \leq s_{mP} \cdot s_0$, so
$$
h(m'P)  =  2 s_{m'P} \cdot s_0  \leq  2 s_{mP} \cdot s_0 = h(mP)
$$
as claimed.  \qed

\textit{Remarks}:
\vspace*{-2ex}
\begin{enumerate}
\item
We could also state the result as:
The \naive\ height of a point is less than or equal to
the \naive\ height of any of its multiples that is not the zero point.
This is a more natural formulation (the first point does not have
to be written as $m'P$\/), but less convenient for our purposes.
\item
In the proof, ``at least the same multiplicity'' can be strengthened
to ``exactly the same multiplicity'' in our characteristic-zero setting.
In general $h(mP)$ may strictly exceed $h(m'P)$ because
$s_{mP} \cap s_0$ may also contain points where $m'P$\/ reduces
to a nontrivial \hbox{$(m/m')$-torsion} point.
\end{enumerate}

The \naive\ height satisfies further inequalities
along the lines of Lemma~\ref{lemma:ht_ineq}, for instance
\be
\label{eq:h6}
h(6P)+h(P) \geq h(2P)+h(3P).
\ee
Lemma~\ref{lemma:ht_ineq} suffices for the proofs of Theorems 1--3
in the genus-zero case,
but inequalities such as~(\ref{eq:h6}) are sometimes needed
to exclude possible configurations with positive~$g$,
as we shall see for $d=24$.
The strongest such inequality we found is:
\begin{lemma}
Let $P$\/ be a point on an elliptic curve over $k(C)$,
and let $m$ be any integer such that $mP \neq 0$.
Then
\label{lemma:ht_mu}
\be
\label{eq:ht_mu}
\sum_{m'|m} \mu(m/m') \, h(m'P) \geq 0.
\ee
\end{lemma}

\textit{Proof}\/:
The left-hand side can be interpreted as
twice the number of points of~$C$, counted with multiplicity,
at which $mP=0$ but $m'P \neq 0$ for each proper factor~$m'$ of~$m$.
\qed

Inequality (\ref{eq:h6}) is the special case $m=6$ of this Lemma.
The sum in (\ref{eq:ht_mu}) may be considered as an analogue of
the formula $\prod_{m'|m} (x^{m'}-1)^{\mu(m/m')}\9$
for the \hbox{$m$-th} cyclotomic polynomial.
We recover Lemma~\ref{lemma:ht_ineq} by summing
the inequality~(\ref{eq:ht_mu}) over all factors of~$m$,
including $m$ itself but not~$1$, to obtain $h(mP) \geq h(P)$,
which is equivalent to Lemma~\ref{lemma:ht_ineq}
by the first Remark above.

\textbf{2.2 Local invariants, and Shioda's inequality.}
To go from the \naive\ to the canonical height we must use
the minimal model of~$E$\/ for the elliptic surface~$\EE$.
We next describe this model, collect some known facts on
the singular fibers of~$\EE$, and give Shioda's lower bound
on the conductor degree.

Whereas a \naive\ height could be defined for any model of~$E$,\footnote{
  Two models may yield different heights $h,h'$, but
  $h'=h+O(1)$ holds for any pair of \naive\ heights
  on the same curve.  It also follows that the property
  $\Ht=h+O(1)$ of the canonical height does not depend
  on the choice of \naive\ height~$h$.
  }
the canonical height requires the N\'eron minimal model.
It is known that there exists a minimal line bundle~$\LL$ on~$C$\/
with the following property: let $D$\/ be a divisor on~$C$\/
such that $O(D) \cong \LL$; then $E$\/ is isomorphic to a curve
with an extended Weierstrass equation (\ref{eq:Weier})
whose coefficients $a_i$ are global sections of~$iD$.
In characteristic zero we can easily obtain~$D$\/ and~$\LL$
by putting~$E$\/ in narrow Weierstrass form
$Y^2 = X^3 + a_4 X + a_6$.  Then $D$\/ is the smallest divisor
such that $(a_4) + 4D \geq 0$ and $(a_6) + 6D \geq 0$.
In other words, we can regard $a_4,a_6$ as global sections
of $\LL^{\otimes 4},\LL^{\otimes 6}$ such that there is no point of~$C$\/
where $a_4$ and $a_6$ vanish to order at least~$4$ and~$6$ respectively.
Once we have $a_i \in \Gamma(\LL^{\otimes i})$, we can regard
the Weierstrass equation (\ref{eq:Weier}) as a surface in
the plane bundle $\LL^{\otimes 2} \oplus \LL^{\otimes 3}$ over~$C$.
If all the roots of the discriminant $\Delta\in\Gamma(\LL^{\otimes 12})$
are distinct then this surface is smooth and is the minimal model of~$E$.
Otherwise it has isolated singularities, which we blow up as many times
as needed (we may follow Tate's algorithm~\cite{Tate_Antwerp})
to obtain the minimal model~$\EE$.  This is a smooth algebraic surface
of arithmetic genus $n = \deg\LL$, equipped with a map to~$C$\/
with generic fiber~$E$\/ and $\omega_{\EE/C} \cong \LL$.
See for instance \cite[pp.149ff.]{BPV}.

We shall need much information about the singular fibers
that can arise for the elliptic fibration $\EE \ra C$.
We extract from Tate's table~\cite[p.46]{Tate_Antwerp}
the following local data for each possible Kodaira type
of a singular fiber~$E_v$: the discriminant degree $d_v$,
the conductor degree $N_v$, and the structure of the group
$E_v/(E_v)\0_0$ of \hbox{multiplicity-$1$} components.
We also list in each case the root lattice~$L_v$
that $E_v$ contributes to the \NeS\ lattice $\NS(\EE)$ of~$\EE$.
In each case, $L_v$ has rank $d_v-N_v$,
and $E_v/(E_v)\0_0 \cong L_v^*/L_v\0$
where $L_v^* \subset L_v\0 \otimes \Q$ is the dual lattice.
The lattice ``$\!A_0$'' that appears for Kodaira types I$_1$ and~II
is the trivial lattice of rank~zero.
For Kodaira type I$^*_\nu$, the group $E_v/(E_v)\0_0$
always has order~$4$, and has exponent~$2$ or~$4$
according as $\nu$ is even or odd.
For positive~$\nu$ of either parity, a fiber of type I$^*_\nu$
has a distinguished \hbox{multiplicity-$1$} component
of order~$2$ in $E_v/(E_v)\0_0$, namely the one closest
to the identity component.  In the $L_v$ picture,
the distinguished component corresponds
to the nontrivial coset of~$D_{4+\nu}$ in~$\Z^{4+\nu}$.
When $\nu=0$ there is no distinguished component:
all three non-identity components of multiplicity~$1$
are equivalent, as are all three nontrivial cosets
due to the triality of~$D_4$.

\vspace*{4ex}

\centerline{
\begin{tabular}{c|cccccccc}
Kodaira type & I$_\nu (\nu>0)$ & II & III & IV &
   \z I$^*_\nu$ & \z IV$^*$ & \z III$^*$ & \z II$^*$
\\ \hline
$d_v$ & $\nu$ & 2 & 3 & 4 & $6+\nu$ & 8 & 9 & 10
\\
$N_v$ &   1   & 2 & 2 & 2 &    2    & 2 & 2 &  2
\\
$E_v/(E_v)\0_0$ & $\Z/\nu\Z$ & $\{0\}$ & $\Z/2\Z$ & $\Z/3\Z$
  &  $\,D_{4+\nu}^*/D_{4+\nu}\0\,$ & $\Z/3\Z$ & $\Z/2\Z$ & $\{0\}$
\\
root lattice & $A_{\nu-1}$ & $A_0$ & $A_1$ & $A_2$ &
   $D_{4+\nu}$ & $E_6$ & $E_7$ & $E_8$
\end{tabular}
}

\vspace*{2ex}

The discriminant and conductor degrees $d,N$\/ of~$\EE$ are sums
of the discriminant and conductor degrees of the singular fibers:
\be
\label{eq:d,N}
12n = d = \sum_v d_v,
\qquad
N = \sum_v N_v.
\ee
Hence $d-N = \sum_v(d_v-N_v) = \sum_v \rk L_v$
is the rank of the subgroup $\oplus_v L_v$ of~$\NS(\EE)$
due to the singular fibers.
Shioda used this to prove \cite[Cor.~2.7 (p.30)]{Sh:Ineq}:

\begin{proposition}
\label{prop:shioda}
Let $E$\/ be a nonconstant elliptic curve over a function field
$K=k(C)$ of genus~$g$, with discriminant and conductor degrees
$d=12n$ and~$N$.  Then
\be
\label{eq:ShIneq}
N \geq 2n + (2-2g) + r,
\ee
where $r$ is the rank of the \MW\ group $E(K)$.
\end{proposition}

\textit{Proof}\/:
Let $T \subseteq \NS(\EE)$ be the subgroup spanned by
$s_0$, the generic fiber, and $\oplus_v L_v$.
Then we have a short exact sequence
(see for instance \cite[Thm.~1.3]{Sh:MW}):
\be
\label{eq:MWseq}
0 \ra T \ra \NS(\EE) \ra E(K) \ra 0,
\ee
where the map $\NS(\EE) \ra E(K)$ is the sum on the generic fiber.
Taking ranks, we find
\be
\label{eq:MWrk}
\rk \NS(\EE) = \rk T + \rk E(K) = 2 + (d-N) + r.
\ee
But $\NS(\EE)$ embeds into $H^{1,1}(\EE,\Z)$,
a group of rank $h^{1,1}(\EE) = 10 n + 2 g$.
Hence $\rk \NS(\EE) \leq 10n+2g$.  Therefore
$$
N \geq (d+2+r) - (10n+2g) = 2n + (2-2g) + r,
$$
as claimed.
\Qed

\textit{Remarks}\/:
\vspace*{-2ex}
\begin{enumerate}
\item
Since $r \geq 0$ it follows that
\be
\label{eq:Szpiro}
N \geq 2n + (2-2g) = (d/6) + \chi
\ee
  for any nonconstant elliptic surface.  This weaker inequality
  is sufficient for most of our purposes, even though we are interested
  in curves with a nontorsion point, for which
  the strict inequality $N > (d/6) + \chi$ holds because $r>0$.
\item
The inequality (\ref{eq:Szpiro}) is now usually known as
  the ``Szpiro inequality''$\!$, but Shioda's paper~\cite{Sh:Ineq}
  predates Szpiro's~\cite{Szpiro} by almost two decades
  (see also \cite[p.114]{Sh:Ineq92}).
  It is by now well-known that (\ref{eq:Szpiro}) can be proved
  by elementary means via Mason's theorem~\cite{Mason}
  (the ABC inequality for function fields).  Can one also give
  an elementary proof of Shioda's inequality,
  or even of its consequence that $r=0$ if $N = (d/6) + \chi$?
\item
The requirement that $E$\/ not be a constant curve is essential.
  There is an analogous statement for constant curves but many details
  must change.  Suppose $E$\/ is such a curve, that is,
  $\EE = C \times E_0$ for some elliptic curve $E_0/k$.
  Then $E(K)$ is not finitely generated, because it contains
  a copy of $E_0(k)$.  Still, $E(K)/E_0(k)$ is finitely generated,
  and identified with the group $\NS(\EE)/T$.
  Again we call the rank of this group~$r$.
  Since $n=d=N=0$ in this setting, we obtain the inequality
  $r + 2 \leq h^{1,1}(C \times E_0) - 2$.  But for a constant curve,
  $h^{1,1}(C \times E_0) = 2g+2$, instead of the $2g$
  that one would expect from the $10n+2g$ formula.  Hence $r \leq 2g$.
  This can also be proved using the identification of $E(K)/E_0(k)$
  with $\End(\Jac(C),E_0)$, an approach that also yields the
  equality condition: clearly $r=2g$ if $g=0$; if $g>0$
  then $r=2g$ if and only if $E_0$ has complex multiplication
  and $\Jac(C)$ is isogenous with $E_0^g$.  See for instance~\cite{E1}.
\item
The hypothesis of characteristic zero, too, is essential here.
  In positive characteristic, one cannot decompose the second
  Betti number~$b_2(\EE)$ as $h^{2,0} + h^{1,1} + h^{0,2}$,
  so one has only the weaker upper bound $b_2(\EE)$ on $\rk(\NS(\EE))$.
  This upper bound exceeds the characteristic-zero bound by
  $2g$ for a constant curve and $2(n+g-1)$ for a nonconstant one.
  For instance, a constant curve $C \times E_0$ has $r \leq 4g$,
  with equality if and only if either $g=0$ or $E_0$ and $\Jac(C)$
  are both supersingular.  In general $\EE$\/ is said to be
  ``supersingular'' if $NS(\EE) \cong \Z^{b_2(\EE)}$;
  such surfaces were studied and used in~\cite{Sh:MW,E1}.
\end{enumerate}

\textbf{2.3 Local height corrections.}
We next list the local height corrections $\lambda_v(mP)$
for each of the Kodaira types.  For convenience we abuse notation
by using $mP$\/ to refer also to the section $s_{mP}$.
\begin{itemize}
\item If $mP$\/ is on the identity component of $E_v$ then
\be
\label{eq:d/6}
\lambda_v(mP) = d_v/6.
\ee
In particular this covers fibers of type~II or~II$^*$.
\item If $E_v$ is of type I$_\nu$ and $P$\/ passes through
component $a \in \Z/\nu\Z$,
let $x=\bar a/\nu$ for any lift~$\bar a$ of~$a$ to~$\Z$; then
\be
\label{eq:nuB}
\lambda_v(mP) = \nu B(mx),
\ee
where $B(\cdot)$ is the second Bernoulli function
 $
 B(z) := \sum_{n=1}^\infty \cos(2\pi n)/(\pi n)^2.
 $
Since $B$\/ is $\Z$-periodic, the choice of $\bar a$ does not matter.
Likewise, since \hbox{$B(z)=B(-z)$} it does not matter that
$a$ cannot be canonically distinguished from $-a$.
We have
\be
\label{eq:B01}
B(z) = z^2 - z + \frac16
\ee
for all $z \in [0,1]$, so in particular $B(0)=1/6$.
Hence $\lambda_v(mP) = \nu/6$
if $mP$\/ passes through the identity component of~$E_v$,
as also asserted by~(\ref{eq:d/6}) in that case.
\item If $E_v$ is of type III, IV, I$^*_0$, III$^*$, or IV$^*$,
and $mP$\/ passes through a non-identity component of~$E_v$,
then $\lambda_v(mP) = 0$.
\item Finally, suppose $E_v$ is of type I$^*_\nu$ ($\nu>0$)
and that $mP$\/ passes through a \hbox{non-identity} component.
If that component is the distinguished one of order~$2$
then $\lambda_v(mP) = \nu/6$.  Otherwise $\lambda_v(mP) = -\nu/12$.
(We could have also allowed $\nu=0$, when there is no distinction
among the three \hbox{non-identity} components, but
$\lambda_v(mP)=\nu/6=-\nu/12=0$ for all of them.)
\end{itemize}



We record two applications of these formulas for future use:

\begin{lemma}
\label{lemma:-n,2n}
Let $E$\/ be an elliptic curve of discriminant degree~$12n$
over a function field~$K$, and $P$\/ any nonzero point of $E(K)$.
Then
\be
\label{eq:-n,2n}
-n \leq  \Ht(P) - h(P) \leq 2n.
\ee
\end{lemma}

\textit{Proof}\/: For each $v$ we have
$-d_v/12 \leq \lambda_v \leq d_v/6$.
Summing over~$v$ yields~(\ref{eq:-n,2n}).
\qed

\begin{lemma}
\label{lemma:mm}
Let $E$\/ be an elliptic curve of discriminant degree~$12n$
over a function field~$K$, and $P$\/ any point of $E(K)$.
If for some integer $m$\/ the multiple $mP$\/ is a nonzero
integral point then $\Ht(mP) \leq 2n/m^2$.
\end{lemma}

\textit{Proof}\/: By our formulas for $\lambda_v$ we have
$\lambda_v(mP) \leq d_v/6$ for all~$v$.  Hence
\be
\label{eq:hmm}
m^2 \Ht(P) = \Ht(mP) = h(mP) + \sum_v \lambda_v(mP)
\leq h(mP) + \sum_v d_v/6.
\ee
But $h(mP)=0$ since $mP$\/ is integral,
and $\sum_v d_v/6 = d/6 = 2n$.
Hence $m^2 \Ht(P) \leq 2n$, and the Lemma follows.
\qed

\textbf{2.4 Reduction to the semistable case.}
Recall that an elliptic curve is said to be {\em semistable}
if all its singular fibers are of type $I_\nu$ for some~$\nu$.
Suppose $E/K$\/ is semistable and $P$\/ is a nontorsion point in~$E(K)$.
We associate to $(E,P)$ an element~$\gamma$ of the abelian group~$\G$
of formal \hbox{$\Z$-linear} combinations of orbits of~$\Q$ under
the infinite dihedral group $D_\infty$ generated by $z \mapsto z+1$
and $z \leftrightarrow 1-z$.  We denote by $[z]$ the generator of~$\G$
corresponding to the orbit of~$z$.  Then $\gamma$ is defined as
a sum of local contributions $\gamma_v \in \G$
that record the types~$\nu(v)$ of the singular fibers~$E_v$
and the component $c_v = a(v) \in \Z/(\nu(v))\Z$ of each fiber
that contains~$P$, as follows:
\be
\label{eq:gammav}
\gamma_v :=
\sum_v \gcd(a(v),\nu(v)) \cdot \left[\frac{a(v)}{\nu(v)}\right] .
\ee
Then each of the height corrections $\Ht(mP)-h(mP)$,
as well as the discriminant degree,
are images of $\gamma$ under homomorphisms $\lala_m,\dd$
from~$\G$ to~$\Q$ or~$\Z$, and the conductor is bounded above
by the image of a homomorphism $\NN : \G \ra \Z$.
We define these homomorphisms on the generators of~$\G$
and extend by linearity.
Suppose $\Q \ni z = a/b$\/ with $b>0$ and $\gcd(a,b)=1$.
Note that $b$\/ is an invariant of the action of~$D_\infty$.
Then we set
\be
\label{eq:l,d,n}
\lala_m([z]) := b \, B_2(mz),
\qquad
\dd([z]) := b,
\qquad
\NN([z]) := 1.
\ee
Then our formulas (\ref{eq:nuB},\ref{eq:d,N}) yield the identities
\be
\label{eq:gamma}
\Ht(mP) = h(mP) + \lala_m(\gamma)  \quad(m=1,2,3,\ldots),
\qquad
12n = d = \dd(\gamma)
\ee
and the estimate
\be
\label{eq:gammaN}
N \leq \NN(\gamma).
\ee
(This last is an upper bound rather than an identity because
each~$v$ contributes~$1$ to~$N$\/ and $\gcd(a(v),\nu(v)) \geq 1$
to~$\NN(\gamma)$.)  It follows that
\be
\label{eq:Nineq}
\NN(\gamma) \geq N \geq (d/6) + (2-2g) + r
\geq \frac16 \dd(\gamma) + 3-2g.
\ee
The second step is Shioda's inequality (Prop.~\ref{prop:shioda}),
and the third step uses the positivity of~$r$, which follows from
our hypothesis that $P$\/ is nontorsion.

To generalize these formulas to curves that may not be semistable,
it might seem that we would have to extend $\G$ with generators
that correspond to Kodaira types other than~$I_\nu$.
But we can associate to any additive fiber~$E_v$ an element of~$\G$
whose images under $\lala_m$ and~$\dd$ coincide with $\lambda_v(mP)$
and~$d_v$, and whose image under~$\NN$ is $\geq N_v$.
(Note that we already did this for multiplicative fibers with
$f = \gcd(a(v),\nu(v)) > 1$, replacing them in effect
by $f$\/ fibers with $a,\nu$ coprime and the same value of~$a/\nu$.)
As in the multiplicative case, this element is positive,
in the sense that it is a nonzero formal linear combination
of elements of $\Q/D_\infty$ with nonnegative coefficients.
Specifically, we have:
\begin{proposition}
\label{prop:G}
Let $E$\/ be an elliptic curve over a function field~$K$\/
of genus~$g$, and $P\in E(K)$ a nontorsion point.
Define for each singular fiber $E_v$ a positive $\gamma_v \in \G$,
depending on $(E_v,c_v)$ as follows:
\begin{itemize}
\item
 If $E_v$ is multiplicative, $\gamma_v$ is defined by~(\ref{eq:gammav}).
\item
 If $c_v$ is the identity component then $\gamma_v := d_v \, [0]$.
\item
 If $c_v$ is a non-identity component of a fiber $E_v$ of type
 III, IV, IV$^*$, or III$^*$ then $\gamma_v$ is respectively
 $$
  [1/2] + [0], \quad
  [1/3] + [0], \quad
  2 \cdot [1/2] + 2 \cdot [0], \quad
  3 \cdot [1/3] + 3 \cdot [0].
 $$
\item
 If $c_v$ is a distinguished component of a fiber $E_v$ of type
 I$^*_\nu$ then
 $$
  \gamma_v := 2 \, [1/2] + (\nu+2) \, [0].
 $$
\item
 If $c_v$ is a non-distinguished, non-identity component
 of a fiber $E_v$ of type I$^*_\nu$ then
 $$
  \gamma_v := (\mu+2) \, [1/2] + 2 \, [0]
 $$
 if $\nu = 2\mu$, and
 $$
  \gamma_v := [1/4] + (\mu+1) \, [1/2] + [0]
 $$
 if $\nu = 2\mu+1$ for some integer~$\mu$.
\end{itemize}
Then:
\nl
$\phantom{Y?}$\
i) $\lambda_v(mP) = \lala_m(\gamma_v)$ for each $m=1,2,3,\ldots$;
\nl
$\phantom{Y?}$\
ii) $d_v = \dd(\gamma_v)$; and
\nl
$\phantom{Y?}$\
iii) $N_v \leq \NN(\gamma_v)$.
\nl
Thus (\ref{eq:gamma},\ref{eq:gammaN},\ref{eq:Nineq}) hold
for $\gamma := \sum_v \gamma_v$.  Equality in~(iii) holds
if and only if $E_v$ is either
a multiplicative fiber with $\gcd(a,\nu)=1$,
a fiber of type III or~IV with $c_v$ a \hbox{non-identity} component,
or a fiber of type~II.
\end{proposition}

[Note that, as was true for the $\lambda_v$ formulas,
the first two formulas in Prop.~\ref{prop:G} overlap
in the case of a multiplicative fiber with $a(v)=0$,
but give the same answer in this case.  Here both prescriptions yield
$\gamma_v = \nu(v) \cdot [0]$ for such~$v$.]

\textit{Proof}\/: The multiplicative case was seen already.
For each of the other Kodaira types, it is straightforward to verify
that $\lambda_v(mP) = \lala_m(\gamma_v)$ for each nonnegative $m$
less than the exponent of the finite group $E_v/(E_v)_0$
(which is at most~$4$), and to check that $d_v = \dd(\gamma_v)$,
and that $N_v \leq \NN(\gamma_v)$,
with strict inequality except in the three cases listed.
We recover (\ref{eq:gamma},\ref{eq:gammaN},\ref{eq:Nineq})
by summing over~$v$.
\Qed

\vspace*{2ex}

{\large\textbf{3.~The values of $\Htmin(0,12n)$ for $n=1,2,3$,
and consecutive integral multiples.}}

For each $n$ we can use the formulas and results above
to obtain a lower bound on $\Htmin(g,12n)$.
When $g=0$ and $n=1,2,3$ we also show that this bound is attained
if and only if $mP$ is integral for $m\leq M=6,8,9$,
and verify that the $(E,P)$ exhibited in Theorem~$n$
satisfy those conditions.

Suppose $E$\/ is an elliptic curve over $\C(T)$
with discriminant degree~$12n$.  Let $P$\/ be a nontorsion
rational point on~$E$, and $\gamma$ the associated element of~$\G$.
{}From $\gamma$ and $\Ht(P)$ we can recover all the \naive\ heights
$h(mP)$ from the first formula in~(\ref{eq:gamma}):
  $
  h(mP) = m^2 \Ht(P) - \lala_m(\gamma).
  $
Given $n$ and an upper bound~$H$\/ on~$\Ht(P)$,
there are only finitely many candidates for the pair $(\gamma,\Ht(P))$:
there are finitely many $\gamma>0$ with $\dd(\gamma)=12n$,
and for each one there are only finitely many possible choices
for $h(P)$ consistent with $h(P) + \lala_1(\gamma) = \Ht(P) \in (0,H]$.
For each candidate $(\gamma,\Ht(P))$ we can check the condition
$m'|m \Ra h(mP) \geq h(m'P) \geq 0$.
Only finitely many $m$ need be checked for each $(\gamma,\Ht(P))$:
by Lemma~\ref{lemma:-n,2n} we know that $h(mP)\geq0$
once $m^2 \Ht(P) \geq n$, and $h(mP) \geq h(m'P)$
for each $m'|m$ once $m^2 \Ht(P) \geq 4n$.
The minimal $\Ht(P)$ among the $(\gamma,\Ht(P))$ that pass these tests
is then our lower bound on~$\Htmin(g,12n)$.
[We could also test the more complicated inequality
of Lemma~\ref{lemma:ht_mu}, which may further improve the bound;
instead we checked that inequality after the fact when necessary.]

We wrote a {\sc gp} program to compute this bound by exhaustive search,
and ran it with $H=2n/M^2$ for $n=1,2,3$.  We chose this upper bound~$H$\/
to ensure that, by Lemma~\ref{lemma:mm}, we would also find all feasible
$(\gamma,\Ht(P))$ such that $h(mP)=0$ for each $m=1,2,3,\ldots,M$.
For $n=1$, we found that the minimum occurs for
\be
\label{eq:min1}
\gamma = [1/5] + [1/3] + [1/2] + 2 \, [0],
\qquad
\Ht(P) = 1/30,
\ee
and is the unique $(\gamma,\Ht(P))$ such that $h(mP)=0$
for each $m \leq 6$.
For $n=2$, we found that the minimum occurs for
\be
\label{eq:min2_not}
\gamma = [1/11] + 2 \, [2/5] + [1/3],
\qquad
\Ht(P) = 4/165;
\ee
but this is not feasible because $h(mP)=0,2,2,2$ for $m=2,4,6,12$,
so inequality~(\ref{eq:h6}) is violated when $m=2$.
Our lower bound on $\Htmin(g,24)$ is thus the next-smallest value,
which occurs for
\be
\label{eq:min2}
\gamma = [1/7] + [2/5] + [1/4] + [1/3] + [1/2] + 3 \, [0],
\qquad
\Ht(P) = 11/420,
\ee
and is the unique $(\gamma,\Ht(P))$ such that $h(mP)=0$
for each $m \leq 8$.

On the other hand, the $(\gamma,\Ht(P))$ pairs
of~(\ref{eq:min1},\ref{eq:min2}) are also those associated with
the curves and points $E,P$\/ exhibited in (\ref{eq:E1},\ref{eq:E2}).
Hence those $E,P$\/ attain our lower bounds $1/30$, $11/420$
on $\Htmin(12)$, $\Htmin(24)$, as well as the upper bounds $6$ and~$8$
on the number of consecutive integral multiples for $n=1$ and $n=2$.
This proves all of Theorems~1 and~2 except for the claims that
every $(E,P)$ attaining those bounds is isomorphic with some
$E_1(q)$ or $E_2(u)$.

For $n=3$, we find that there is a unique $(\gamma,\Ht(P))$
such that $h(mP)=0$ for each $m \leq 9$, namely
\be
\label{eq:min3}
\gamma = [1/8] + [3/7] + [1/5] + [1/4] + 2 \, [1/3] + [1/2] + 4 \, [0],
\quad
\Ht(P) = 23/840.
\ee
Again these are the $\gamma$ and $\Ht(P)$ for the $(E,P)$
exhibited in the Introduction (formula~(\ref{eq:E3})).
But we do not claim that $\Htmin(36)=23/840$:
Lemma~\ref{lemma:ht_mu} eliminates the second-smallest pair
$$
(\gamma,\Ht(P)) = ([1/13]+[3/8]+[3/7]+[1/5]+[1/3], \  229/10920)
$$
(which violates the inequality (\ref{eq:h6})
in the same way that (\ref{eq:min2_not}) did),
but not several other possibilities with $\Ht(P) < 23/840$.
We next list all these possibilities, in order of increasing $\Ht(P)$:
\be
\label{eq:min3_not}
\begin{array}{c|r}
\gamma & \Ht(P) \phantom{0000} \\ \hline
[1/13] + [3/11] + [3/8] + 2\,[1/2]              & 23/1144 \approx .02010 \\ {}
[1/13] + [3/8] + [2/7] + [1/4] + 2\,[1/2]       & 17/728  \approx .02335 \\ {}
[1/11] + [4/9] + [2/7] + [1/4] + [1/3] + 2\,[0] & 65/2772 \approx .02345 \\ {}
[1/12] + [3/11] + [3/8] + 2\,[1/2] + [0]        &  7/264  \approx .02652 \\ {}
[1/11] + [3/7] + 2\,[1/5] + [1/4] + 2\,[1/2]    & 41/1540 \approx .02662
\end{array}
\ee
(For comparison, $229/10920 \approx .02097$ and $23/840 \approx .02738$.)
We have $\dd(\gamma) \leq 7$ for each entry
in the table~(\ref{eq:min3_not});
therefore by Prop.~\ref{prop:shioda} none of them can occur
for an elliptic curve over~$\PP^1$.
(Even the weaker inequality (\ref{eq:Szpiro}) would suffice here;
either of those inequalities also excludes (\ref{eq:min2_not})
for $n=2$, and would thus be enough to obtain $\Htmin(0,24)$,
but the determination of $\Htmin(24)$ required a further argument.)
Thus $\Htmin(0,36) = 23/840$, proving Theorem~3 except for the claim
that every $(E,P)$ satisfying conditions (a) and~(b) is of the form
$E_3(A)$ for some~$A$.

\textbf{Acknowledgements.}
I thank J.~Silverman and T.~Shioda for helpful correspondence
concerning their papers and related issues,
and M.~Watkins for carefully reading a draft of this paper.
This work was made possible in part by funding from
the Packard Foundation and the National Science Foundation.

\end{document}